\documentclass[12pt,reqno]{amsart}

\usepackage{amssymb}

\newcommand{\sbs}{\subset}
\newcommand{\seq}{\subseteq}
\newcommand{\stm}{\setminus}
\newcommand{\est}{\varnothing}

\newcommand{\del}{\delta}
\newcommand{\sig}{\sigma}
\newcommand{\lam}{\lambda}
\renewcommand{\phi}{\varphi}

\newcommand{\N}{{\mathbb N}}
\newcommand{\R}{{\mathbb R}}
\newcommand{\Z}{{\mathbb Z}}

\newcommand{\cI}{{\mathcal I}}
\newcommand{\cL}{{\mathcal L}}
\newcommand{\cX}{{\mathcal X}}

\newcommand{\prt}{\partial}

\newcommand{\longc}{,\ldots,}

\newcommand{\longp}{+\dotsb+}
\newcommand{\longge}{\ge\dotsb\ge}
\newcommand{\longcup}{\cup\dotsb\cup}

\newcommand{\cl}[1]{\operatorname{cl}(#1)}
\newcommand{\scl}[1]{\operatorname{cl^*}(#1)}
\renewcommand{\int}[1]{\operatorname{int}(#1)}
\newcommand{\sint}[1]{\operatorname{int^*}(#1)}

\newcommand{\refl}[1]{\ref{l:#1}}
\newcommand{\reft}[1]{\ref{t:#1}}
\newcommand{\refp}[1]{\ref{p:#1}}
\newcommand{\refs}[1]{\ref{s:#1}}
\newcommand{\refb}[1]{\cite{b:#1}}
\newcommand{\refe}[1]{\eqref{e:#1}}

\newtheorem{lemma}{Lemma}
\newtheorem{theorem}{Theorem}
\newtheorem{primelemma}{Lemma}

\newcounter{step}
\newcommand{\step}%
  {\smallskip\noindent\refstepcounter{step}{\bf Step~\arabic{step}.\ }}

\title[Minimizing the projection sums]%
  {Minimizing the sum of projections \\ of a finite set}

\author{Vsevolod F. Lev}
\email{seva@math.haifa.ac.il}
\address{Department of Mathematics, The University of Haifa at Oranim,
  Tivon 36006, Israel}

\author{Misha Rudnev}
\email{M.Rudnev@bristol.ac.uk}
\address{Department of Mathematics, University of Bristol, Bristol BS8~1TW,
  United Kingdom}

\keywords{Isoperimetric problem, Loomis-Whitney inequality, projections}
\subjclass[2010]{Primary: 52A38; secondary: 52A40, 05A20}

\begin{document}
\baselineskip = 16pt

\begin{abstract}
Consider the projections of a finite set $A\sbs\R^n$ onto the coordinate
hyperplanes. How small can the sum of the sizes of these projections be,
given the size of $A$? In a different form, this problem has been studied
earlier in the context of edge-isoperimetric inequalities on graphs, and it
is can be derived from the known results that there is a linear order on the
set of $n$-tuples with non-negative integer coordinates, such that the sum in
question is minimised for the initial segments with respect to this order. We
present a new, self-contained and constructive proof, enabling us to obtain a
stability result and establish algebraic properties of the smallest possible
projection sum. We also solve the problem of minimising the sum of the sizes
of the one-dimensional projections.
\end{abstract}

\maketitle

\section{Preliminaries}\label{s:prelim}

Given an integer $n\ge 1$, for each $i\in[1,n]$ denote by $\pi_i$ the
orthogonal projection of the vector space $\R^n$ onto the coordinate
hyperplane $\{(x_1\longc x_n)\in\R^n\colon x_i=0\}$. For a finite set
$A\sbs\R^n$, as a simple consequence of the Loomis-Whitney inequality (see,
for instance, \cite{b:lw,b:cgg,b:gmr}), we have
  $$ \prod_{i=1}^n |\pi_i(A)| \ge |A|^{n-1}; $$
combining this estimate with the arithmetic mean-geometric mean inequality yields
\begin{equation}\label{e:LW+AGM}
  \sum_{i=1}^n |\pi_i(A)| \ge n\, |A|^{1-1/n}.
\end{equation}

The Loomis-Whitney inequality is known to be sharp, turning into an exact
equality when the set under consideration is an axes-aligned rectangular
parallelepiped (and, in the convex set settings, only in this case, as it
follows from the argument of \refb{lw}). In contrast, the estimate
\refe{LW+AGM} is not sharp; say, it shows that for a three-dimensional,
five-point set, the sum of the projection sizes is at least as large as
$3\cdot 5^{2/3}\approx 8.77$, while it is not difficult to see that, in fact,
this sum cannot be smaller than $10$. This leads naturally to the following
question: exactly how small can the sum in the left-hand side of
\refe{LW+AGM} be for a finite set $A\sbs\R^n$ of given size? Loosely
speaking, we want to know how much the points of an $n$-dimensional set of
given size can hide behind each other.

The answer to a tightly related question is due to Bollob\'as and Leader
\cite[Theorem~15]{b:bl}, where it was cast as the edge-isoperimetric problem
for the $n$-dimensional grid graph; see also Harper \cite[Theorem~7.1]{b:h}
and historical comments \cite[page 142]{b:h}, as well as the references
contained therein. Our goal here is to give a direct, independent, and
self-contained solution of a discrete-geometric flavour, avoiding references
to graph theory and making the underling rearrangement procedure maximally
algorithmic. This enables us to prove a stability result showing that, in
certain cases, the set with the smallest sum of the projection sizes is,
essentially, unique. Furthermore, we establish some algebraic properties of
the smallest possible value of this sum as a function of the size of the set
being projected. Finally, in the Appendix we discuss and solve a similar
problem for the one-dimensional projections.

\section{Summary of results}\label{s:intro}

If $|A|=K^n$ with an integer $K\ge 1$, then \refe{LW+AGM} implies
  $$ \sum_{i=1}^n|\pi_i(A)| \ge nK^{n-1}, $$
which is attained for the discrete $n$-dimensional cube $A=[0,K-1]^n$. The
situation where $|A|$ is \emph{not} a perfect $n$th power is much subtler and
requires some preparation to discuss.

Denote by $\N_0$ the set of all non-negative integers. Following \refb{bl},
we define the \emph{cube order} on $\N_0^n$ by saying that $(x_1\longc x_n)$
precedes $(y_1\longc y_n)$ if there exists $l_0\in\N_0$ and $j\in[1,n]$ such
that $\{i\in[1,n]\colon x_i=l\}=\{i\in[1,n]\colon y_i=l\}$ for each $l>l_0$,
and also $\{i\in[j+1,n]\colon x_i=l_0\}=\{i\in[j+1,n]\colon y_i=l_0\}$, while
$x_j<y_j=l_0$.
% KEEPING THE LARGE COORDINATES SMALL, LESS CARE ABOUT SMALL COORDINATES
For integer $m\ge 0$, by $\cI_n(m)$ we denote the length-$m$ initial segment
of $\N_0^n$ with respect to the cube order; thus, for instance,
$\cI_n(0)=\est$, $\cI_1(m)=[0,m-1]$,
  $$ \cI_2(10)
       = \{(0,0),(1,0),(0,1),(1,1),(2,0),(2,1),(0,2),(1,2),(2,2),(3,0)\}, $$
and
\begin{align*}
  \cI_3(17) = \{
     &(0,0,0),(1,0,0),(0,1,0),(1,1,0),
                           (0,0,1),(1,0,1),(0,1,1),(1,1,1), \\
     &(2,0,0),(2,1,0),(2,0,1),(2,1,1),(0,2,0),(1,2,0),(0,2,1),
                                                       (1,2,1),(2,2,0) \}.
\end{align*}
Speaking about \emph{initial segments} we will always mean finite initial
segments of $\N_0^n$ with respect to the cube order, with the value of $n$
determined by the context.

We notice that the cube order is quite similar, but not identical to the
order introduced in \cite[Section~7.1.1]{b:h}; in fact, the latter order is
defined on the whole grid $\Z^n$.

We say that the initial segment $I_1$ is \emph{shorter} than the initial
segment $I_2$ if $|I_1|<|I_2|$; equivalently, if $I_1\sbs I_2$.

Informally, the initial segments fill in $\N_0^n$ cube-wise: once the cube
$[0,K-1]^n$ has been filled in, for some integer $K\ge 1$, the $n$ faces
  $$ 0\le x_1\longc x_i\le K,\ x_{i+1}=K,\ 0\le x_{i+2}\longc x_n\le K-1,
                                                      \quad i\in[0,n-1] $$
are completed one by one to get a covering of the next cube $[0,K]^n$, etc.
If $m=(K+1)^iK^{n-i}$ with some $i\in[0,n-1]$, then
\begin{equation}\label{e:rectpar}
  \cI_n(m) = \{ (x_1\longc x_n)\in\N_0^n\colon 0\le x_1\longc x_i\le K,
                                          \ 0\le x_{i+1}\longc x_n\le K-1 \}
\end{equation}
is an axes-aligned rectangular parallelepiped; in this case we say that the
segment $\cI_n(m)$ is \emph{closed} (the intuition behind this term will be
clear from the next section). The \emph{edges} of a closed initial segment
are its orthogonal projections onto the coordinate axes; thus, for instance,
the two-element initial segment has one edge of size $2$, all other edges
being of size $1$.

In general, for any integer $m\ge 1$ there are uniquely defined integers
$K\ge 1$ and $i\in[0,n-1]$ such that
\begin{equation}\label{e:defK}
  K^n\le m<(K+1)^n
\end{equation}
and, indeed,
\begin{equation}\label{e:defi}
  (K+1)^iK^{n-i} \le m<(K+1)^{i+1}K^{n-i-1};
\end{equation}
writing then
\begin{equation}\label{e:defR}
 m=(K+1)^iK^{n-i}+R,\quad 0\le R<(K+1)^iK^{n-1-i},
\end{equation}
the initial segment $\cI_n(m)$ is the disjoint union of the closed initial
segment $\cI_n((K+1)^iK^{n-i})$, which is the parallelepiped in the
right-hand side of \refe{rectpar}, and a translate of the $(n-1)$-dimensional
initial segment $\cI_{n-1}(R)$, contained in the hyperplane $x_{i+1}=K$.

We can now state our first principal result.
\begin{theorem}\label{t:main}
For every integer $n\ge 1$ and every finite set $A\sbs\mathbb\R^n$, letting
$m:=|A|$, we have
\begin{equation}\label{e:sig}
  \sum_{i=1}^n |\pi_i(A)| \ge \sum_{i=1}^n |\pi_i(\cI_n(m))|.
\end{equation}
\end{theorem}

In Section~\refs{pmain} we prove Theorem~\reft{main} in the particular case
where $A\sbs\N_0^n$; the general case then follows readily by observing that
if $A\seq S^n$ with a finite set $S\seq\R$, then for any injective mapping
$\phi\colon S\to\N_0$, writing $B:=\phi^{\otimes n}(A)\sbs \N_0^n$ we have
$|B|=|A|$ and $|\pi_i(A)|=|\pi_i(B)|$ for each $i\in[1,n]$. Indeed, this
observation shows that the estimate of Theorem~\reft{main} remains valid when
$A\seq S^n$ with a set $S$ of any nature, not necessarily contained in the
set of real numbers (although in this case the projections $\pi_i$ must be
redefined appropriately).

We denote the left-hand side of \refe{sig} by $\sig_n(A)$ and, with some
abuse of notation, its right-hand side by $\sig_n(m)$; that is, $\sig_n(A)$
is the sum of the sizes of the projections of the finite set $A\sbs\R^n$ onto
the coordinate hyperplanes, and $\sig_n(m)$ is this sum in the special case
where the set in question is the length-$m$ initial segment. Thus, for
instance,
\begin{equation}\label{e:initcond}
  \sig_n(0)=0,\ \sig_1(m)=1\ \text{if}\ m>0,
\end{equation}
and, as one can easily verify,
\begin{equation}\label{e:sig2}
  \sig_2(m) = \begin{cases}
                    2K+1\ &\text{if}\ K^2<m\le K(K+1), \\
                    2K+2\ &\text{if}\ K(K+1)<m\le (K+1)^2
                 \end{cases}
\end{equation}
for any integer $K\ge 0$. Also, it follows from the explanation above that
for $K,i$, and $R$ defined by \refe{defK}--\refe{defR}, we have
\begin{align}\label{e:recurs}
  \sig_n(m) &= \sig_n((K+1)^iK^{n-i}) + \sig_{n-1}(R) \notag \\
            &= (nK+n-i)(K+1)^{i-1}K^{n-i-1} + \sig_{n-1}(R);
\end{align}
along with \refe{initcond}, this relation gives a recursive, completely
algebraic definition of the quantities $\sig_n(m)$.

We now address the corresponding stability problem.

We say that a finite set $A\sbs\R^n$ is a \emph{minimiser} if its projection
sum $\sig_n(A)$ is smallest possible among all sets in $\R^n$ of size $|A|$.
Thus, Theorem~\reft{main} says that every initial segment of $\N_0^n$ is a
minimiser, but it is not true in general that any minimiser is an initial
segment, or even is ``similar'' to an initial segment in some reasonable
sense (see Section~\refs{pmainu} for a rigorous definition of similarity).
Say, for integer $K\ge C\ge 1$, the set
$A:=[0,K-C]\times[0,K+C]\subset\N_0^2$ is not an initial segment, while
$|A|=(K+1)^2-C^2$ and therefore $\sig_2(A)=2K+2=\sig_2(|A|)$
(cf.~\refe{sig2}), showing in view of Theorem~\reft{main} that $A$ is a
minimiser. In Section~\refs{pmainu} we prove, however, that every
\emph{closed} initial segment is (up to similarity) a unique minimiser.
\begin{theorem}\label{t:mainu}
Suppose that $n\ge 1$ is an integer. If $m=(K+1)^iK^{n-i}$ with integers
 $K\ge 1$ and $i\in[0,n-1]$, then every minimiser in $\R^n$ of size $m$ is a
Cartesian product of $i$ real sets of size $K+1$, and $n-i$ real sets of size
$K$.
\end{theorem}

The following lemma, proved in Section~\refs{psub}, is an important
ingredient of the proof of Theorem~\reft{main}.
\begin{lemma}\label{l:sub}
Let $n\ge 1$ be an integer.
\begin{itemize}
\item[i)] Suppose that $I_1,I_2,J_1,J_2\sbs\N_0^n$ are initial segments
    such that $|I_1|+|I_2|=|J_1|+|J_2|$, $|J_1|\le|I_1|\le|I_2|\le|J_2|$,
    and  $J_2$ is closed. Then
      $$ \sig_n(J_1)+\sig_n(J_2) \le \sig_n(I_1)+\sig_n(I_2). $$
\item[ii)] If $I,I_1,I_2\sbs\N_0^n$ are non-empty initial segments such
    that $|I|=|I_1|+|I_2|$, then
      $$ \sig_n(I) < \sig_n(I_1)+\sig_n(I_2). $$
\item[iii)] If $n\ge 2$ and $I_{n-1}\sbs\N_0^{n-1}$, $I_n\sbs\N_0^n$ are
    non-empty initial segments such that $|I_{n-1}|=|I_n|$, then
    $$ \sig_n(I_n) > \sig_{n-1}(I_{n-1}). $$
\end{itemize}
\end{lemma}

An ``algebraic restatement'' of Lemma~\refl{sub} may be of interest.
\begin{primelemma}\label{l:sub-algebraic}
Let $n\ge 1$ be an integer.
\begin{itemize}
\item[i)] Suppose that $l_1\le m_1\le m_2\le l_2$ are non-negative
    integers such that $m_1+m_2=l_1+l_2$. If $\cI_n(l_2)$ is closed, then
      $$ \sig_n(l_1) + \sig_n(l_2)\le \sig_n(m_1)+\sig_n(m_2). $$
\item[ii)] If $m_1,m_2\ge 1$ are integers, then
      $$ \sig_n(m_1+m_2) < \sig_n(m_1)+\sig_n(m_2). $$
\item[iii)] If $n\ge 2$ and $m\ge 1$ is an integer, then
    $$ \sig_n(m) > \sig_{n-1}(m). $$
\end{itemize}
\end{primelemma}

We remark that subadditivity established by Lemma~\refl{sub} ii) and
Lemma~\refl{sub-algebraic} ii) can be viewed as a combinatorial analogue of
the well-known physical fact that merging two spherical droplets into one
reduces the total surface area.

Our last result establishes yet another interesting algebraic property of the
functions $\sig_n$.
\begin{theorem}\label{t:restate}
For any integers $n,s\ge 1$ and $m_1\longc m_s\ge 0$, we have
\begin{equation}\label{e:additivity}
  \sig_n(m_1\longp m_s) \le \sig_{n-1}(m_1)\longp \sig_{n-1}(m_s)
                                                 + \max \{ m_1\longc m_s \}.
\end{equation}
\end{theorem}
In Section~\refs{prestate} we derive Theorem~\reft{restate} from
Theorem~\reft{main} and, indeed, show that, somewhat unexpectedly, the two
theorems are equivalent in the sense that each of them follows easily from
the other one.

In the next section we introduce important notation and terminology used
throughout. Having finished with this, we prove Lemma~\refl{sub} in
Section~\refs{psub}, and Theorems~\reft{main} and~\reft{mainu} in
Sections~\refs{pmain} and~\refs{pmainu}, respectively. The equivalence of the
former theorem and Theorem~\reft{restate} is established in
Section~\refs{prestate}. Finally, in the Appendix we develop a similar theory
for the one-dimensional projections, and in particular show that their sum is
also minimised when the set under consideration is an initial segment.

The proofs are purely combinatorial, based on point rearrangements.

\section{Notation and terminology}\label{s:notation}

For integer $1\le i\le n$, we denote the coordinate vectors in $\R^n$ by
$e_1\longc e_n$, and write
  $$ \cX_i := \operatorname{Sp} \{e_i\} $$
for the coordinate axes, and
  $$ \cL_i := \operatorname{Sp}\{e_1\longc e_{i-1},e_{i+1}\longc e_n\} $$
for the coordinate hyperplanes. The axis $\cX_n$ will be referred to as
\emph{vertical}, and the corresponding hyperplane $\cL_n$ and its translates,
as well as the projection $\pi_n$, as \emph{horizontal}.

The reader may find helpful to think of points $x\in\Z^n$ as unit cubes, with
the base vertex at $x$, and visualize sets $A\seq\Z^n$ as built of such
cubes.

The intersections of a set $A\seq\Z^n$ with horizontal hyperplanes will be
called the \emph{slabs} of $A$. If $n=1$, then the slabs are
zero-dimensional; hence, either empty, or singletons. Notice that closed
initial segments in $\N_0^n$ are stable under permuting non-empty slabs.

Recall that an initial segment of $\N_0^n$ is \emph{closed} if it is an
axes-aligned rectangular parallelepiped whose edges differ by at most $1$,
and for all $1\le i<j\le n$, the edge along $\cX_j$ is not longer than the
edge along $\cX_i$ (cf.~\refe{rectpar}).

Given an initial segment $I\sbs\N_0^n$, we define its (strict)
\emph{interior} to be the longest closed initial segment (strictly) contained
in $I$, and we denote the interior and the strict interior of $I$ by
$\int{I}$ and $\sint{I}$, respectively. Similarly we define the (strict)
\emph{closure} of $I$ to be the shortest closed initial segment (strictly)
containing $I$, and denote the closure and the strict closure of $I$ by
$\cl{I}$ and $\scl{I}$, respectively. The boundary of $I$ is defined by
 $\prt I:=I\stm\int{I}$, and its strict boundary by $\prt^*I:=I\stm\sint{I}$;
thus, the boundary is empty if and only if $I$ is closed, while the strict
boundary is always nonempty whenever $I\ne\est$. Boundaries can be
treated either as $(n-1)$-dimensional sets embedded in $\N_0^n$, or as
initial segments in $\N_0^{n-1}$.

We remark that for an initial segment $I\ne\est$, any of the three conditions
$\scl{I}=\cl{I}$, $\sint{I}=\int{I}$, and $\prt^*I=\prt I$ is equivalent to
$I$ not being closed.

As a version of \refe{recurs}, for any initial segment $I\sbs\N_0^n$ with
$|I|>1$, we have
\begin{equation}\label{e:idt}
  \sig_n(I) = \sig_n(\sint{I}) + \sig_{n-1}(\prt^*I).
\end{equation}
(If $|I|=1$, then $\sint{I}=\est$ and the left-hand side of \refe{idt}
exceeds by $1$ its right-hand side.) This basic, but important identity
allows us to argue inductively in the forthcoming proofs. It becomes evident
upon observing that the strict boundary $\prt^*I$ is an $(n-1)$-dimensional
set, attached to and not larger than a face of $\sint{I}$ which we visualize
as a rectangular parallelepiped; hence $\prt^*I$  does not contribute to the
projection along the axis, normal to the face under consideration.

Notation-wise, we will occasionally use \refe{idt} in the form
\begin{equation*}\label{e:idtm}
  \sig_n(|I|) = \sig_n(|\sint{I}|) + \sig_{n-1}(|\prt^*I|),
\end{equation*}
the equivalence following from the fact that $\sint{I}$ and $\prt^*I$ are
initial segments (in $\N_0^n$ and $\N_0^{n-1}$, respectively).

\section{Proof of Lemma~\refl{sub}}\label{s:psub}

We use induction by $n$. The base case $n=1$ is easy to verify, and we
proceed assuming that $n\ge 2$. We first prove iii), then i), and, finally,
ii).

Addressing iii), we use (the inner) induction by the common size $m$ of the
initial segments $I_{n-1}$ and $I_n$. For $m=1$ the estimate follows from
  $$ \sig_n(I_n) = n = \sig_{n-1}(I_{n-1})+1. $$
For $m\ge 2$, by \refe{idt} and the induction hypothesis, we have
  $$ \sig_n(I_n) = \sig_n(|\sint{I_n}|) + \sig_{n-1}(|\prt^* I_n|)
                    > \sig_{n-1}(|\sint{I_n}|) + \sig_{n-1}(|\prt^* I_n|). $$
Applying ii) inductively, we see that the right-hand side is larger than
  $$ \sig_{n-1}(|\sint{I_n}|+|\prt^* I_n|) = \sig_{n-1}(|I_n|)
                            = \sig_{n-1}(|I_{n-1}|) = \sig_{n-1}(I_{n-1}). $$
This completes the proof of iii), and we now turn to i).

If $J_1=\est$, then the assertion follows from ii), and we thus assume that
$J_1\ne\est$, implying $I_1\ne\est$. Our plan is to shorten $I_1$ and
lengthen $I_2$, while keeping the sum $|I_1|+|I_2|$ intact, to get, after a
number of iterations, to the situation where $I_2=J_2$. Formally, we act as
follows.

Let
\begin{equation}\label{e:del}
  \del := \min \{ |\prt^*I_1|, |\scl{I_2}\stm I_2| \}
\end{equation}
(the number of elements to be transferred from $I_1$ to $I_2$), and define
$I_1'$ and $I_2'$ to be the initial segments of $\N_0^n$ of sizes
$|I_1'|=|I_1|-\del$ and $|I_2'|=|I_2|+\del$. Notice that $\del>0$, and since
$J_2$ is closed, if $I_2\ne J_2$, then we have $\scl{I_2}\seq J_2$, implying
$\del\le|J_2\stm I_2|$. Consequently, $|I_2'|=|I_2|+\del\le|J_2|$, whence
$|I_1'|=|J_1|+|J_2|-|I_2'|\ge|J_1|$.

We now prove that
\begin{equation}\label{e:mdashsigma}
  \sig_n(I_1') + \sig_n(I_2') \le \sig_n(I_1) + \sig_n(I_2).
\end{equation}
If $I_2'\ne J_2$ (equivalently, if $I_1'\ne J_1$), then we iterate the
procedure, until eventually we replace the initial segments $I_1$ and $I_2$
with $J_1$ and $J_2$, respectively, and the assertion will then follow from
\refe{mdashsigma}. Thus, to complete the proof of i) it remains to establish
\refe{mdashsigma}. To this end, we have to distinguish two cases.

Suppose first that $|\prt^*I_1|>|\scl{I_2}\stm I_2|$, so that $I_2$ is not
closed in view of $|I_1|\le|I_2|$; consequently, $|\prt I_2|=|\prt^* I_2|>0$
and $\int{I_2}=\sint{I_2}$. In this case we have $\del=|\scl{I_2}|-|I_2|$
whence $|I_2'|=|\scl{I_2}|$ and therefore $I_2'=\scl{I_2}$ and
$\sint{I_2'}=\sint{I_2}$; also, $|I_1'|>|I_1|-|\prt^* I_1|=|\sint{I_1}|$,
implying $\sint{I_1'}=\sint{I_1}$. As a result, using \refe{idt}, we get
\begin{align}\label{e:tri}
  \sig_n(I_1) + \sig_n(I_2)
    &= \sig_n(\sint{I_1}) + \sig_{n-1}(\prt^* I_1) + \sig_n(\sint{I_2})
                                         \notag + \sig_{n-1}(\prt^* I_2) \\
    &= \sig_n(\sint{I_1'}) +  \sig_n(\sint{I_2'})
         + \big( \sig_{n-1}(\prt^* I_1) + \sig_{n-1}(\prt^* I_2) \big).
\end{align}
We now notice that $\prt^*I_2'$ is a closed initial segment in $\N_0^{n-1}$,
and that
  $$ |\prt^*I_1'|+|\prt^*I_2'| = (|\prt^*I_1|-\del) + (|\prt^*I_2|+\del)
                                               = |\prt^*I_1|+|\prt^*I_2|. $$
Also,
  $$ |\prt^*I_2'| = |I_2'|-|\sint{I_2'}| = |\scl{I_2}|-|\sint{I_2}|
                                      > |I_2|-|\sint{I_2}| = |\prt^*I_2|. $$
Therefore, an inductive application of i) in dimension $n-1$ yields
  $$ \sig_{n-1}(\prt^* I_1) + \sig_{n-1}(\prt^* I_2)
               \ge \sig_{n-1}(\prt^* I_1') + \sig_{n-1}(\prt^* I_2'). $$
Combining this with \refe{tri}, and using \refe{idt} once again, we obtain
\begin{align*}
  \sig_n(I_1) + \sig_n(I_2)
    &\ge \sig_n(\sint{I_1'}) + \sig_n(\sint{I_2'})
              + \sig_{n-1}(\prt^* I_1') + \sig_{n-1}(\prt^* I_2') \\
    &=   \sig_n(I_1') + \sig_n(I_2'),
\end{align*}
which is the desired estimate~\refe{mdashsigma}.

Now suppose that $|\prt^*I_1|\le|\scl{I_2}\stm I_2|$. In this case
$\del=|\prt^*I_1|$, $I_1'=\sint{I_1}$, $\sint{I_2'}=\int{I_2}$, and
\begin{equation}\label{e:boundaries}
  |\prt^*I_2'| = |I_2'\stm\int{I_2}|
                         = |\prt I_2| + \del = |\prt^* I_1| + |\prt I_2|.
\end{equation}
There are two further sub-cases.

If $I_2$ is closed, then \refe{boundaries} gives $|\prt^*I_2'|=|\prt^* I_1|$;
as a result, using \refe{idt} we get
\begin{align}
  \sig_n(I_1) + \sig_n(I_2)
    &\ge \big( \sig_n(\sint{I_1})
               + \sig_{n-1}(\prt^*I_1) \big) + \sig_n(I_2) \label{e:odin} \\
    &= \sig_n(I_1')
           + \big( \sig_{n-1}(\prt^*I_1) + \sig_n(\int{I_2}) \big) \notag \\
    &= \sig_n(I_1') + \big( \sig_{n-1}(\prt^*I_2')
                                       + \sig_n(\sint{I_2'}) \big) \notag \\
    &= \sig_n(I_1') + \sig_n(I_2'), \notag
\end{align}
which is \refe{mdashsigma}. (The inequality in \refe{odin} is strict if and
only if $I_1$ is a singleton; this fact will be used in the forthcoming proof
of ii).)

If $I_2$ is \emph{not} closed, then $\prt I_2=\prt^*I_2$ and
$\sint{I_2}=\int{I_2}=\sint{I_2'}$. Recalling~\refe{boundaries}, in this case
we apply ii) inductively in dimension $n-1$ to get
  $$ \sig_{n-1}(\prt^* I_1) + \sig_{n-1} (\prt^* I_2)
                                                > \sig_{n-1}(\prt^* I_2') $$
whence, by \refe{idt},
\begin{align}
  \sig_n(I_1) + \sig_n(I_2)
     &\ge \sig_n(\sint{I_1}) + \big( \sig_n(\sint{I_2})
            + \sig_{n-1}(\prt^*I_1) + \sig_{n-1}(\prt^* I_2) \big) \notag \\
     &>   \sig_n(I_1') + \big( \sig_n(\sint{I_2'})
                            + \sig_{n-1}(\prt^* I_2') \big) \label{e:dva} \\
     &=   \sig_n(I_1')+\sig_n(I_2'). \notag
\end{align}
This establishes~\refe{mdashsigma}, and therefore~i).

Finally, we prove ii). Without loss of generality, we assume $|I_1|\le|I_2|$.
If
\begin{equation}\label{e:scl2large}
  |\scl{I_2}\stm I_2| \le |I_1|,
\end{equation}
then we re-use the above argument for i) with $J_2:=\scl{I_2}$ and
$J_1:=\cI_n(|I_1|+|I_2|-|J_2|)$, defining $\del$ by \refe{del} and then
letting $I_1':=\cI_n(|I_1|-\del)$ and $I_2':=\cI_n(|I_2|+\del)$, to have the
estimate \refe{mdashsigma}. We notice that if $|I_1|>1$, then $\del\le|\prt^*
I_1|<|I_1|$, implying $|I_1'|\ge 1$; moreover, if $|I_1|=1$, then the
inequality in \refe{mdashsigma} is strict, as it follows from \refe{dva} and
the remark following \refe{odin}. (This reflects the geometrically obvious
fact that if $I_1$ consists of one single point, then moving this point to
$I_2$ reduces the total sum of the hyperplane projections by at least $1$.)

Continuing in this way, we find initial segments $I_1''\seq I_2''$ satisfying
$|I_1''|+|I_2''|=|I_1|+|I_2|$ and
\begin{equation}\label{e:primeprime}
  \sig_n(I_1'')+\sig_n(I_2'') \le \sig_n(I_1)+\sig_n(I_2)
\end{equation}
such that either $I_1''=\est$, $I_2''=I$, and \refe{primeprime} holds
actually as a strict inequality, or $I_1''\ne\est$ and $|\scl{I_2''}\stm
I_2''|>|I_1''|$, cf.~\refe{scl2large}. In the former case ii) follows
readily. In the latter case, recalling that $I=\cI_n(|I_1''|+|I_2''|)$, we
have $\sint{I}=\int{I_2''}$, whence
  $$ |\prt^*I| = (|I_1''|+|I_2''|) - |\int{I_2''}| = |\prt I_2''|+|I_1''|; $$
consequently, using \refe{idt}, and then applying~ii) inductively,
\begin{align}
  \sig_n(I)
    &=   \sig_n(\sint{I}) + \sig_{n-1}(\prt^* I) \notag \\
    &=   \sig_n(\int{I_2''}) + \sig_{n-1}(|\prt I_2''|+|I_1''|) \notag \\
    &\le \sig_n(\int{I_2''})
               + \sig_{n-1}(\prt I_2'') + \sig_{n-1}(I_1''). \label{e:temp1}
\end{align}
However, as a version of \refe{idt} (essentially equivalent
to~\refe{recurs}), we have
\begin{equation}\label{e:temp2}
  \sig_n(\int{I_2''}) + \sig_{n-1}(\prt I_2'') = \sig_n(I_2''),
\end{equation}
and by iii),
\begin{equation}\label{e:temp3}
  \sig_{n-1}(I_1'') < \sig_n(I_1'').
\end{equation}
Combining \refe{temp1}--\refe{temp3} and \refe{primeprime}, we get
  $$ \sig_n(I) < \sig_n(I_1'') + \sig_n(I_2'')
                                           \le \sig_n(I_1) + \sig_n(I_2). $$
This completes the proof of ii).

\section{Proof of Theorem~\reft{main}}\label{s:pmain}

As explained in the introduction, it suffices to show that for any finite set
$A\sbs\N_0^n$, writing $m:=|A|$, we have $\sig_n(A)\ge\sig_n(m)$. The proof
goes by induction on $n$, the base case $n=1$ being trivial as the
zero-dimensional projection of any nonempty set has by convention, size $1$.
The assertion is readily verified for $m\in\{0,1\}$, too. Suppose thus that
$\min\{n,m\}\ge 2$.

Our strategy is to start out with any minimiser $A\sbs\N_0^n$ and modify it,
in a finite number of rearrangements not increasing the projection sum, to
get the initial segment $\cI_n(|A|)$. We achieve this in several steps, some
of which may need to be iterated, as explained below.

\step\label{p:step1}%
Let $H$ (for ``height'') denote this number of non-empty slabs of $A$. We
permute the slabs so that the number of elements of $A$ in any higher slab
does not exceed the number of elements in a lower slab; that is, letting
$A^{(k)}:=A\cap(ke_n+\cL_n)$, we have $|A^{(0)}|\longge|A^{(H-1)}|>0$ and
$|A^{(k)}|=0$ for $k\notin[0,H-1]$. This does not affect the projection sum
$\sig_n(A)$ as each non-horizontal projection of $A$ is the disjoint union of
the corresponding projections of the slabs:
\begin{equation}\label{e:HZ1}
  \pi_i(A) = \bigcup_{k\ge 0} \pi_i(A^{(k)}), \quad i\in[1,n-1].
\end{equation}
To simplify the notation, we keep denoting by $A$ the set under
consideration. Clearly, the number of non-empty slabs of $A$ remains equal to
$H$.

We now replace each slab $A^{(k)}$ with the initial segment
$\cI_{n-1}(|A^{(k)}|)$, without enlarging the projection sum. (It is readily
seen that the horizontal projection does not increase, and the sum of the
side projections does not increase in view of~\refe{HZ1} and by the induction
hypothesis.) We use the same notation $A$ for the new set obtained in this
way, but from now on we assume that each slab of $A$ is an
$(n-1)$-dimensional initial segment. Since the sequence $(|A^{(k)}|)_{k\ge
0}$ is non-increasing, this implies $A^{(k+1)}\seq e_n+A^{(k)}$ for each
$k\ge 0$.

Let $K$ denote the largest edge size of the closed $(n-1)$-dimensional
initial segment $\cl{A^{(0)}}$; that is, the size of the projection of
$A^{(0)}$ onto the coordinate axis $\cX_1$. If $K<H$, then we swap the
coordinate axes $\cX_1$ and $\cX_n$ so that the number of non-empty slabs
decreases to $K$, and repeat the whole procedure.

We keep permuting the slabs and swapping the axes until $A$ gets rearranged
so that the number $H$ of non-empty slabs does not exceed the largest edge
size $K$ of the closure $\cl{A^{(0)}}$ of the lowest slab, and each slab is
an $(n-1)$-dimensional initial segment.

\step\label{p:step2}%
A repeated application of this step will ensure that all non-empty slabs of
$A$, with the possible exception of the highest slab $A^{(H-1)}$, have the
same interior. Assuming this is not the case, there are integers
$k\in[1,H-2]$ with $|A^{(k)}|<|\int{A^{(0)}}|$. Let $k$ be the smallest such
integer. If, indeed, we had $|A^{(k)}|+|A^{(H-1)}|\le|\int{A^{(0)}}|$, then
we would be able to remove $A^{(H-1)}$ from $A$ and replace $A^{(k)}$ with
the (appropriate vertical translate of the) initial segment
$\cI_{n-1}(|A^{(k)}|+|A^{(H-1)}|)$, without changing the horizontal
projection $\pi_n(A)$. By Lemma~\refl{sub} ii), this would result in the
strict decrease of the sum of the non-horizonal projections, contradicting
the assumption that $A$ is a minimiser. Thus, we have
$|A^{(k)}|+|A^{(H-1)}|>|\int{A^{(0)}}|$, and we replace the slab $A^{(k)}$
with $\int{A^{(0)}}$, and the upper slab $A^{(H-1)}$ with the initial segment
$\cI_{n-1}(|A^{(H-1)}|+|A^{(k)}|-|\int{A^{(0)}}|)$; by Lemma~\refl{sub} i),
applied with $I_1=A^{(H-1)}$, $I_2=A^{(k)}$, and $J_2=\int{A^{(0)}}$, this
does not increase the sum of non-horizontal projections of $A$, and it is
clear that the horizontal projection $\pi_n(A)$ remains unchanged.

We emphasize that Step~\refp{step2} affects neither the number $H$ of the
slabs of $A$, nor the lower slab $A^{(0)}$, and that if this step ever gets
applied, then the resulting set satisfies $|A^{(0)}|>|A^{(H-1)}|$.

Repeating Step~\refp{step2}, we ensure that all non-empty slabs of $A$,
excepting perhaps $A^{(H-1)}$, have their interiors identical to that of
$A^{(0)}$, which we assume to hold from now on.

\step\label{p:step3}%
Recall that by $K$ we denote the largest edge size of the closed
$(n-1)$-dimensional segment $\cl{A^{(0)}}$. As a result of the rearrangements
of Step~\refp{step1}, we have $K\ge H$, and the present Step~\refp{step3} is
to be repeated as long as the strict inequality $K>H$ holds, or until $A$
gets rearranged as in initial segment.

If $A$ can be cast as an initial segment just by relabelling the coordinate
axes, this will complete the proof; this scenario will be referred to as a
\emph{trivial exit}. Otherwise, we are going to cut from $A$ a ``vertical
slab'' resting on the strict boundary $\prt^*A^{(0)}$, and place it as a new
horizontal slab (as a result of which $H$ will grow by $1$). Given that only
``side'' projections of both slabs contribute to the projection sum
$\sig_n(A)$, it will not be affected by this rearrangement. Formally, we need
to consider two cases.

The first case is $|A^{(H-1)}|\ge|\sint{A^{(0)}}|$ (covering, in particular,
the situation where $H=1$). In this case each slab of $A$ contains a vertical
translate of the strict interior $\sint{A^{(0)}}$, and we define $A'$ to be
the union of all these $H$ translates, and let $A'':=A\stm A'$. Notice that
$A''$ is a non-empty subset of a hyperplane parallel to one of the
(non-horizontal) coordinate hyperplanes. Considering $A''$ as an
$(n-1)$-dimensional set, we have
\begin{equation}\label{e:onemore}
  \sig_n(A) = \sig_n(A') + \sig_{n-1}(A'').
\end{equation}
Observe that if we replace the set $A''$ with the initial segment
$\cI_{n-1}(|A''|)$, by the induction assumption and the assumption that $A$
is a minimiser, this will not affect the projection sum $\sig_{n-1}(A'')$.

If $H=K-1$, then $A'$ is the union of $K-1$ vertical translates of the closed
initial segment $\sint{A^{(0)}}$, which is the $n$-dimensional axes-aligned
rectangular parallelepiped with the edge sizes $K-1$ and (possibly) $K$. The
set $A''$ is strictly smaller than the face of $A'$ it is attached to, for
otherwise the original set $A$ would be a closed initial segment, up to
relabelling of the coordinate axes, and we would have the trivial exit
scenario. Hence $A''$ can be replaced with the initial segment
$\cI_{n-1}(|A''|)$ and re-attached to the appropriate face of $A'$ to get a
set which, up to a coordinate axes relabelling, is an initial segment,
completing the proof (for the present subcase $H=K-1$). Observe, for the
forthcoming proof of Theorem \reft{mainu}, that the fact that $A''$ is
strictly smaller than the face of $A'$ it is attached to, precludes the
output initial segment $\cI_n(|A|)$ from being closed.

Assuming now that $H<K-1$, let $K_1\longge K_{n-1}$ be the edge sizes of
$\sint{A^{(0)}}$, so that $K\ge K_1$ and $K_{n-1}\ge K-1$. We have then
\begin{equation}\label{e:A''issmall}
  |A''| \le |\prt^* A^{(0)}|\cdot  H < K_1\dotsb K_{n-2}\cdot (K-1)
                      \le K_1\dotsb K_{n-2}\cdot K_{n-1} = |\sint{A^{(0)}}|.
\end{equation}
It follows that we can detach $A''$ from $A$, replace it with the initial
segment $\cI_{n-1}(|A''|)$, and re-attach as the upper slab, thus increasing
$H$ by $1$ and changing the lower slab of $A$ from the original $A^{(0)}$ to
its strict interior $\sint{A^{(0)}}$. Notice that, as a consequence of
\refe{A''issmall}, the resulting set has fewer elements in its upper slab
than in the lower slab.

For the second case $|A^{(H-1)}|<|\sint{A^{(0)}}|$, we define $A'$ to be the
union of $H-1$ (rather than $H$ as in the first case) vertical translates of
the set $\sint{A^{(0)}}$, with the coordinate $\cX_n$ ranging from $0$ to
$H-2$, and we let $A'':=A\stm(A'\cup A^{(H-1)})$. Thus $A'$ is an
axes-aligned rectangular parallelepiped, attached to two faces of which are
the $(n-1)$-dimensional sets $A''$ and $A^{(H-1)}$; hence,
  $$ \sig_n(A) = \sig_n(A') + \sig_{n-1}(A'') + \sig_{n-1}(A^{(H-1)}). $$
We now rearrange $A''$ as an $(n-1)$-dimensional initial segment, and then
detach it from $A$ and re-attach as either the upper, or the
second-from-the-top slab, to retain the non-increasing order of slab sizes.
As above, this rearrangement makes $H$ larger by $1$, and changes the lower
slab of $A$ from the original $A^{(0)}$ to its strict interior
$\sint{A^{(0)}}$. Also, keeping denoting the notation for the edges of
$\sint{A^{(0)}}$, similarly to~\refe{A''issmall} we have
\begin{equation}\label{e:A''issmalll}
  |A''| \le |\prt^* A^{(0)}|\cdot (H-1) \le K_1\dotsb K_{n-2}\cdot (H-1)
                       < K_1\dotsb K_{n-2}\cdot K_{n-1} = |\sint{A^{(0)}}|;
\end{equation}
thus, as above, the size of the horizontal projection of $A$ becomes
$|\sint{A^{(0)}}|$.

Finally (just for the second-from-the-top slab) we invoke the rearrangement
of Step~\refp{step2} to ensure that all, but the upper slab of $A$ have the
same interior; hence, are actually identical closed initial segments since
$A^{(0)}$ is a closed initial segment.

Observe that, unless we have achieved our goal of rearranging $A$ as an
initial segment (as in the trivial exit scenario or the case where
$|A^{(H-1)}|\ge|\sint{A^{(0)}}|$ and $H=K-1$), the procedure introduced in
Step~\refp{step3} results in $H$ growing by $1$, with the strict inequality
$|A^{(H-1)}|<|A^{(0)}|$ for the rearranged $A$, and with $A^{(0)}$ being a
closed $(n-1)$-dimensional initial segment. In addition, the parameter $K$ is
either unchanged, or decreases by $1$, the latter happening if and only if
the new ``base slab'' $A^{(0)}$ is an $(n-1)$-dimensional cube. Therefore, if
the new parameters satisfy $K<H$ (that is, $H=K+1$), then the set $A$ got
rearranged into a cube with the edge size $K$, with an $(n-1)$-dimensional
initial segment attached to its upper face as boundary; that is, into a
(non-closed) $n$-dimensional initial segment.

We have shown that, applying Step~\refp{step3} repeatedly, we will either
rearrange $A$ as in initial segment, or arrive in the special situation where
$K=H$, dealt with at Step~\refp{step4} below.

\step\label{p:step4}%
For this last step of our algorithm we assume that $K=H$ where, we recall,
$H$ is the number of non-empty slabs of $A$, and $K$ is the largest size of a
projection of the lower slab $A^{(0)}$ onto a non-vertical coordinate axis.
This step is not to be iterated; it is applied at most once and after
completing it, $A$ will be rearranged as an initial segment.

If $A$ can be rearranged as an initial segment by merely relabelling the
coordinate axes, then the algorithm stops and proof is completed; as in
Step~\refp{step3}, this situation will be referred to as the \emph{trivial
exit}. Assuming that we have not exited trivially, as in the above Step 3, we
consider two cases.

The first one is a straightforward modification of the corresponding case of
Step~\refp{step3}, with $\sint{A^{(0)}}$ replaced by $\int{A^{(0)}}$ (and the
equality $H=K-1$ replaced by $H=K$). Namely, if
$|A^{(H-1)}|\ge|\int{A^{(0)}}|$, then we define $A'$ to be the union of $H=K$
vertical translates of $\int{A^{(0)}}$, one on top of the other, and let
$A''=A\setminus A'$. The set $A'$ is an $n$-dimensional axes-aligned
rectangular parallelepiped with the maximum edge size $K$ and minimum edge
size at least $K-1$, and $A''$ is attached to a maximal-sized face of $A'$.
Hence, having $A''$ replaced by the $(n-1)$-dimensional initial segment
$\cI_{n-1}(|A''|)$, the set $A$ can be cast as the initial segment by
relabelling the coordinate axes.

Preparing the ground for the proof of Theorem~\reft{mainu} in the next
section, we notice that $A''$ is non-empty, and is strictly smaller than the
face of $A'$ it is attached to, for otherwise $A$ can be rearranged as an
initial segment by relabelling the coordinate axes, which would lead to the
trivial exit scenario. This precludes the output set $\cI_n(|A|)$ from being
closed.

Moving on to the second case, for the rest of Step~\refp{step4} we assume
that $|A^{(H-1)}|<|\int{A^{(0)}}|$. The set $A\stm A^{(H-1)}$ consists of
$H-1=K-1$ slabs, each one being an $(n-1)$-dimensional initial segment with
the same interior $\int{A^{(0)}}$. All projections onto the non-vertical axes
of the closed $(n-1)$-dimensional initial segment $\int{A^{(0)}}$ have size
$K$ or $K-1$; therefore, the stack of $K-1$ vertical translates of
$\int{A^{(0)}}$, which we denote $A'$, is a closed $n$-dimensional initial
segment.  Furthermore, the set $A'':=A\stm(A'\cup A^{(H-1)})$ is nonempty,
for otherwise $A$ could be rearranged as the initial segment by relabelling
the coordinate axes, which is ruled out by trivial exit scenario.

Thus, attached to the upper horizontal face of $A'$ is the slab $A^{(H-1)}$,
and to some of its ``vertical'' faces --- the ``vertical slab'' $A''$ (which,
by the construction, is strictly smaller than the face of $A'$ it is attached
to). Our plan is to replace $A''$ by the same-sized $(n-1)$-dimensional
initial segment, and then apply Lemma~\refl{sub} to transfer elements from
$A^{(H-1)}$ and $A''$ to $A'$, to augment this latter set to its strict
closure $\scl{A'}$.

Since the horizontal projections of $A^{(H-1)}$ and $A''$ are disjoint, we
have
  $$ \sig_n(A) = \sig_n(A') + \sig_{n-1}(A'') + \sig_{n-1}(A^{(H-1)}). $$
By the induction assumption, replacing $A''$ with the initial segment
$\cI_{n-1}(|A''|)$ does not increase the summand $\sig_{n-1}(A'')$. As usual,
we do not change the notation, but assume below that $A''$ is an
$(n-1)$-dimensional initial segment.

Let $I_1$ and $I_2$ be the smallest and the largest among the initial
segments $A''$ and $A^{(H-1)}$, respectively, and let $J_2:=\scl{A'}\stm A'$;
thus, $J_2$ is the face to be added to the closed initial segment $A'$ in
order to obtain the ``next'' closed initial segment $\scl{A'}$. Notice that,
by the virtue of the cube order, $|J_2|$ is the size of the \emph{largest}
face of $A'$, whence $|J_2|>|I_2|$. Also notice that $I_1$ and $I_2$ are
non-empty as so are $A''$ and $A^{(H-1)}$. If we had $|J_2|\ge|I_1|+|I_2|$
then, applying Lemma~\refl{sub}~ii), we could have replaced $A''$ and
$A^{(H-1)}$ with one single $(n-1)$-dimensional initial segment of size
$|A''|+|A^{(H-1)}|$ attached to the largest face of $A$, decreasing the sum
of the projections of $A$; this would contradict the assumptions that $A$ is
a minimiser. Therefore we have $|I_1|\le|I_2|\le|J_2|<|I_1|+|I_2|$, and we
set $J_1:=\cI_{n-1}(|I_1|+|I_2|-|J_2|)$ and replace $A''$ and $A^{(H-1)}$
with the initial segments $J_2$ attached to the appropriate face of $A'$ to
convert it to $\scl{A'}$, and $J_1$ attached as a boundary to $\scl{A'}$
(which can be done in view of $|J_1|\le|J_2|$). This rearranges $A$ as an
initial segment. Observe that for this last exit scenario, the output set
$\cI_n(|A|)$ cannot be closed in view of the estimate $0<|J_1|<|J_2|$
resulting from
  $$ |J_1| = |I_1|+|I_2|-|J_2| \le 2|I_2|-|J_2| < |J_2|. $$

\section{Proof of Theorem~\reft{mainu}}\label{s:pmainu}

Applying the argument presented after the statement of Theorem~\reft{main} in
Section~\refs{intro}, we assume without loss of generality that
$A\sbs\N_0^n$.

We define \emph{similarities} to be bijective transformations of the set
$\Z^n$ involving (finitely many) permutations of the horizontal hyperplanes,
axes relabellings, and compositions thereof. Thus, two sets in $\Z^n$ are
\emph{similar} if they can be obtained from each other by a finite series of
permutations of the slabs and relabellings of the coordinate axes. For $n=1$,
any two sets of the same size are similar, proving the assertion in this
case. For $n\ge 1$, if $A_1,A_2\sbs\Z^n$ are similar, finite sets, then
$\sig_n(A_1)=\sig_n(A_2)$.

Our argument uses induction by $n$ and is based on a careful examination of
the proof of Theorem~\reft{main} in the previous section; in fact, it has
been prepared by the observations made there, and particularly at the key
Steps~\refp{step3} and~\refp{step4}. Clearly, it suffices to show that if
$A\sbs\N_0^n$ is a minimiser such that $\cI_n(|A|)$ is closed, then all
rearrangements made in the course of the proof are, in fact, similarities;
that is, involve only permuting the slabs and relabelling the coordinates.

Suppose thus that $A\seq\N_0^n$ is a minimiser with $\cI_n(|A|)$ closed.
Inspecting the proof of Theorem~\reft{main}, we make the two following
observations.
\begin{itemize}
\item[(i)] When Step~\refp{step3} is applied with a non-trivial exit, the
    output set has strictly fewer elements in its upper slab than in the
    lower one; in particular, this set cannot be converted into a closed
    initial segment by relabelling the axes.
\item[(ii)] The only way that Steps~\refp{step3} and~\refp{step4} can
    yield a closed initial segment is that they are exited trivially; in
    particular, the input set must be an axes-aligned rectangular
    parallelepiped.
\end{itemize}
It follows that after completing Steps~\refp{step1} and~\refp{step2}, the set
$A$ may have only required an axes relabelling to get transformed into a
closed initial segment. In fact, no application of Step 2 would have been
possible either, for any such application results in a set with its upper
slab strictly smaller than the lower one.

Thus, rearranging $A$ to $\cI_n(|A|)$ has only required Step~\refp{step1}
followed, possibly, by an axes relabelling. We recall that Step~\refp{step1}
consists of a number of iterations of the procedure that involves permuting
slabs, replacing each slab with the equal-sized $(n-1)$-dimensional initial
segment, and swapping the axes.

Consider the last iteration of Step~\refp{step1}; specifically, the middle
part of the iteration replacing each slab with the same-sized
$(n-1)$-dimensional initial segment. Let $A'$ and $A''$ denote the
corresponding input and output sets. Thus, $A''$ is an axes-aligned
rectangular parallelepiped, with all of its edges differing in size by $1$ at
most. It follows that the common size of all slabs of $A''$ is the
cardinality of a closed $(n-1)$-dimensional initial segment, and we invoke
the induction hypothesis to conclude that on the last iteration, replacing
each slab of $A'$ with an $(n-1)$-dimensional initial segment is induced by
an $(n-1)$-dimensional similarity transformation.

On the other hand, we note that the horizontal projection of $A'$ is the
union of its slabs, viewed as $(n-1)$-dimensional sets (which have the same
common size, as so do the slabs of $A''$). Since we are working with
minimisers, this implies that the size of this union is equal to the size of
each individual slab of $A'$, and therefore all the slabs of $A'$ are
actually identical. As a result, the same similarity transformation that
converts, say, the lower slab of $A'$ into a closed $(n-1)$-dimensional
initial segment, will also work for all other slabs of $A'$ converting them
into (identical) closed initial segments. Extending this transformation to
act as an identity on the last coordinate, we obtain an $n$-dimensional
similarity transformation that replaces all slabs of $A'$ with the
$(n-1)$-dimensional initial segments.

We conclude that rearranging $A'$ into $A''$, and hence the whole last
iteration of Step~\refp{step1}, can be achieved using a similarity
transformation. Making our way backwards, the same is true for all the
preceding iterations. Consequently, the whole Step~\refp{step1} acted on $A$
as a similarity, and the assertion follows.

\section{Proof of the equivalence of Theorems~\reft{main} and~\reft{restate}}%
  \label{s:prestate}

Given integers $n,s\ge 1$ and $m_1\longc m_s\ge 0$, consider the set
$A\sbs\N_0^n$ with $s$ non-empty slabs such that for every $k\in[1,s]$, the
$k$th slab is the $(n-1)$-dimensional initial segment of length $m_k$. Since
the side projections of all these slabs are pairwise disjoint, while the
horizontal projections are all contained in the largest of them, we have
  $$ \sig_n(A)
          = \sig_{n-1}(m_1)\longp\sig_{n-1}(m_s) + \max\{m_1\longc m_s\}. $$
However, the left-hand side is at least as large as
$\sig_n(|A|)=\sig_n(m_1\longp m_s)$ by Theorem~\reft{main}. Hence
Theorem~\reft{main} implies Theorem~\reft{restate}.

Conversely, assuming Theorem~\reft{restate}, one can prove
Theorem~\reft{main} by induction on $n$, as follows. Given a finite set
$A\sbs\N_0^n$, consider the slab decomposition
$A=A^{(1)}\longcup A^{(s)}$, with the $A^{(k)}$ listing all non-empty slabs
of $A$. For each $k\in[1,s]$, let $m_k:=|A^{(k)}|$. Trivially, we have
$|\pi_n(A)|\ge\max\{m_1\longc m_s\}$. Also, disjointness of the side
projections yields
  $$ \sum_{i=1}^{n-1} |\pi_i(A)|
                        = \sum_{k=1}^s \sum_{i=1}^{n-1} |\pi_i(A^{(k)})|. $$
By the induction hypothesis, the double sum in the right-hand side is at
least as large as $\sig_{n-1}(m_1)\longp \sig_{n-1}(m_s)$. Therefore, by
Theorem~\reft{restate},
\begin{multline*}
  \sum_{i=1}^n |\pi_i(A)|
    \ge \sig_{n-1}(m_1)\longp \sig_{n-1}(m_s) + \max \{ m_1\longc m_s \} \\
                                    \ge \sig_n(m_1\longp m_s) = \sig_n(|A|),
\end{multline*}
as claimed by Theorem~\reft{main}.

\appendix
\section*{Appendix: One-dimensional projections}

It would be interesting to extend our results onto $k$-dimensional
projections for all integers $k\in[1,n-2]$. Below we consider the case $k=1$,
establishing the analogs of Theorems~\reft{main},~\reft{mainu},
and~\reft{restate} for the one-dimensional projections. In particular, we
show that the sum of the sizes of these projections is also minimised on the
initial segments of $\N_0^n$ with respect to the cube order, and prove a
corresponding stability result.

For each $i\in[1,n]$, denote by $\rho_i$ the orthogonal projections of $\R^n$
onto the coordinate axis $\cX_i$, and given a finite set $A\sbs\R^n$, let
  $$ \lam_n(A) := \sum_{i=1}^n |\rho_i(A)|; $$
also, for integer $m\ge 0$ let $\lam_n(m):=\lam_n(\cI_n(m))$. Thus, for
instance, $\lam_1(m)=m$, $\lam_2(m)=\sig_2(m)$, and if $|A|=K^n$ with an
integer $K\ge 1$, then by the arithmetic-geometric mean inequality,
  $$ \lam_n(A) \ge n \left(\prod_{i=1}^n |\rho_i(A)|\right)^{1/n}
                                                    \ge n |A|^{1/n} = nK, $$
with equality attained for the discrete cube $A=[0,K-1]^n$.

The key to understanding the quantity $\lam_n$ is the equality
\begin{equation}\label{e:HZ19}
 \lam_n(m+1) = \begin{cases}
        \lam_n(m) + 1\ &\text{if $\cI_n(m)$ is closed}, \\
        \lam_n(m)    \ &\text{otherwise}.
                   \end{cases}
\end{equation}
An immediate corollary is that if $K^n\le m\le(K+1)^n$, with a positive
integer $K$, and $i\in[0,n]$ is the smallest integer such that $m\le
(K+1)^iK^{n-i}$, then
  $$ \lam_n(m) = \lam_n((K+1)^iK^{n-i}) = nK+i, $$
cf.~\refe{recurs}.

As an analog of Theorems~\reft{main} and~\reft{mainu}, we now have
\begin{theorem}\label{t:k=1}
For every integer $n\ge 1$ and every finite set $A\sbs\mathbb\R^n$, letting
$m:=|A|$, we have $\lam_n(A)\ge\lam_n(m)$. Moreover, if $m=(K+1)^iK^{n-i}$
with integers $K\ge 1$ and $i\in[0,n]$, then equality is attained if and only
if $A$ is the Cartesian product of $i$ real sets of size $K+1$, and $n-i$
real sets of size $K$.
\end{theorem}

Although it is possible to prove Theorem~\reft{k=1} modifying the proofs of
Theorems~\reft{main} and~\reft{mainu} to apply in our present settings,
somewhat surprisingly, one can get away with a much easier, non-inductive
argument.

\begin{proof}[Proof of Theorem~\reft{k=1}]
The case where $m=0$ is trivial, and we assume that $m>0$; that is, $A$ is
non-empty. For each $j\in[1,n]$, let $m_j:=|\rho_j(m)|$; thus,
 $m\le m_1\dotsb m_n$ and $\lam_n(A)=m_1\longp m_n$. If the largest of the
numbers $m_j$ exceeds the smallest of them by at least $2$, then we decrease
by $1$ the largest, and simultaneously increase by $1$ the smallest; clearly,
this operation does not affect the sum of the numbers, and their product will
only get larger. Iterating, we will eventually find $n$ positive integers,
the largest of them exceeding the smallest one by at most $1$, so that their
product is at least $m$, and their sum is $\lam_n(A)$. Denoting by $I$ the
closed initial segment whose edges are determined by these resulting
integers, we then have $m\le |I|$ and $\lam_n(A)=\lam_n(I)$. The former
relation gives $\cI_n(m)\seq I$, and then the latter yields
$\lam_n(A)=\lam_n(I)\ge\lam_n(\cI_n(m))=\lam_n(m)$. This proves the first
assertion of the theorem.

For the second assertion, assume that $\lam_n(A)=\lam_n(m)=nK+i$. We also
assume without loss of generality that $i\le n-1$ (if $i=n$, then we can
replace $K$ with $K+1$). Analyzing the argument above, we conclude that if
the smallest among the projections $|\rho_j(A)|$ differed from the largest by
at least $2$, then the size of the initial segment $I$ would satisfy the
strict inequality $m<|I|$, implying $\lam_n(A)=\lam_n(I)>\lam_n(m)$ in view
of \refe{HZ19}, since $\cI_n(m)$ is closed, a contradiction. It follows that
the largest of the projections $|\rho_j(A)|$ differs from the smallest one by
at most $1$. Let $L$ denote the smallest of these projections, and let
$k\in[0,n-1]$ be the number of indices $j\in[1,n]$ with $|\rho_j(A)|=L+1$ (so
that there are $n-k$ those indices $j\in[1,n]$ with $|\rho_j(A)|=L$). From
$nK+i=\lam_n(A)=(L+1)k+L(n-k)=nL+k$ and $k,i\in[0,n-1]$ we then conclude that
$L=K$ and $k=i$. Thus, $A$ is contained in the Cartesian product of $i$ sets of
size $K+1$ and $n-i$ sets of size $K$, and it is in fact equal to this
product since $|A|=(K+1)^iK^{n-i}$.
\end{proof}

In conclusion, we remark that the one-dimensional analog of
Theorem~\reft{restate} is the estimate
  $$ \lam_n(sm) \le \lam_{n-1}(m) + s,\quad m,s\ge 0. $$
The interested reader will easily verify that this estimate follows from the
first assertion of Theorem~\reft{k=1} and, in fact, is equivalent to it.

\section*{Acknowledgement}
We thank Ben Barber and Larry Harper who have kindly helped us to rectify our
ignorance about the edge-isoperimetric inequalities.

\bigskip

\end{document}